\magnification 1200
\vsize 9.0 truein

\input amssym.def
\overfullrule0pt
\hyphenation{pre-print}
\nopagenumbers
\headline={\ifnum\pageno=1 \hfill \else\hss{\tenrm--\folio--}\hss \fi}

\newcount\sectionnumber
\newcount\equationnumber
\newcount\thnumber
\newcount\refnumber

\def\ifundefined#1{\expandafter\ifx\csname#1\endcsname\relax}
\def\assignnumber#1#2{%
	\ifundefined{#1}\relax\else\message{#1 already defined}\fi
	\expandafter\xdef\csname#1\endcsname
 {\the\sectionnumber.\the#2}}%
\def\secassignnumber#1#2{%
  \ifundefined{#1}\relax\else\message{#1 already defined}\fi
  \expandafter\xdef\csname#1\endcsname{\the#2}}%
%
%
\def\secname#1{\relax
  \global\advance\sectionnumber by 1
  \secassignnumber{S#1}\sectionnumber
  \csname S#1\endcsname}
\def\Sec#1 #2 {\vskip0pt plus.1\vsize\penalty-250\vskip0pt plus-.1\vsize
  \bigbreak\bigskip
  \equationnumber0\thnumber0
  \noindent{\bf \secname{#1}. #2}\par
  \nobreak\smallskip\noindent}
\def\sectag#1{\ifundefined{S#1}\message{S#1 undefined}{\sl #1}%
  \else\csname S#1\endcsname\fi}
%
%
\def\eq#1{\relax
  \global\advance\equationnumber by 1
  \assignnumber{EN#1}\equationnumber
  {\rm (\csname EN#1\endcsname)}}
\def\eqtag#1{\ifundefined{EN#1}\message{EN#1 undefined}{\sl (#1)}%
  \else(\csname EN#1\endcsname)\fi}
%
%
\def\thname#1{\relax
  \global\advance\thnumber by 1
  \assignnumber{TH#1}\thnumber
  \csname TH#1\endcsname}
\def\thtag#1{\ifundefined{TH#1}\message{TH#1 undefined}{\sl #1}%
  \else\csname TH#1\endcsname\fi}
\def\Assumption#1 {\bLP{\bf Assumption \thname{#1}}\quad}
\def\Cor#1 {\bLP{\bf Corollary \thname{#1}}\quad}
\def\Def#1 {\bLP{\bf Definition \thname{#1}}\quad}
\def\Example#1 {\bLP{\bf Example \thname{#1}}\quad}
\def\Lemma#1 {\bLP{\bf Lemma \thname{#1}}\quad}
\def\Prop#1 {\bLP{\bf Proposition \thname{#1}}\quad}
\def\Remark#1 {\bLP{\bf Remark \thname{#1}}\quad}
\def\Theor#1 {\bLP{\bf Theorem \thname{#1}}\quad}
\def\Proof{\LP{\bf Proof}\quad}
%
%
\def\refitem#1 #2\par{\ifundefined{REF#1}
\global\advance\refnumber by1%
\expandafter\xdef\csname REF#1\endcsname{\the\refnumber}%
\else\item{\ref{#1}}#2\sLP\fi}

\def\ref#1{\ifundefined{REF#1}\message{REF#1 is undefined}\else
  [\csname REF#1\endcsname]\fi}
\def\reff#1#2{\ifundefined{REF#1}\message{REF#1 is undefined}\else
  [\csname REF#1\endcsname, #2]\fi}
\def\Ref{\vskip0pt plus.1\vsize\penalty-250\vskip0pt plus-.1\vsize
  \bigbreak\bigskip\leftline{\bf References}\nobreak\smallskip
  \frenchspacing}
%
%
\let\sPP=\smallbreak
\let\mPP=\medbreak
\let\bPP=\bigbreak
\def\LP{\par\noindent}
\def\sLP{\smallbreak\noindent}
\def\mLP{\medbreak\noindent}
\def\bLP{\bigbreak\noindent}

\def\halmos{\hbox{\vrule height0.31cm width0.01cm\vbox{\hrule height
 0.01cm width0.3cm \vskip0.29cm \hrule height 0.01cm width0.3cm}\vrule
 height0.31cm width 0.01cm}}
\def\hhalmos{{\unskip\nobreak\hfil\penalty50
	\quad\vadjust{}\nobreak\hfil\halmos
	\parfillskip=0pt\finalhyphendemerits=0\par}}

\def\al{\alpha}
\def\be{\beta}
\def\ga{\gamma}
\def\de{\delta}
\def\la{\lambda}
\def\si{\sigma}

\def\De{\Delta}

\def\CC{{\Bbb C}}
\def\RR{{\Bbb R}}
\def\ZZ{{\Bbb Z}}

\def\FSD{{\cal D}}

\def\gog{{\goth g}}
\def\goh{{\goth h}}

\def\iy{\infty}
\def\id{{\rm id}}
\def\Zplus{\ZZ_{\ge0}}
\def\half{{\scriptstyle{1\over2}}}
\def\thalf{{\textstyle{1\over2}}}
\let\lan=\langle
\let\ran=\rangle
\let\ten=\otimes
\let\union\cup
\def\tfrac#1#2{{\textstyle{#1\over#2}}}

\def\wt{{\rm wt}}
\let\dt\cdot
\def\Tr{{\rm Tr}}

\font\titlefont=cmr10 at 17.28truept
\font\titletexti=cmmi10 at 17.28truept
\font\titletextrm=cmr10 at 17.28truept
\font\authorfont=cmti10 at 14truept
\font\addressfont=cmr10 at 10truept
\font\ttaddressfont=cmtt10 at 10truept
 
\refitem A-A-R

\refitem Bab

\refitem B-B-B

\refitem Erd1

\refitem E-S

\refitem E-V1

\refitem E-V2

\refitem E-V3

\refitem E-V4

\refitem Fe1

\refitem Fe2

\refitem G-N

\refitem Gou

\refitem K-S

\refitem Koo84

\refitem Kra

\refitem Rio

\refitem W

{\titlefont\textfont1\titletexti \textfont0\titletextrm\LP
QDYBE: some explicit formulas for exchange matrix and
\mLP
related objects in case of $sl(2)$, $q=1$}
\bLP
{\authorfont Tom H. Koornwinder and Nabila Touhami}
\bLP
Version of July 14, 2000

\bLP
{\bf Abstract}
\LP
This mainly tutorial paper is intended as a somewhat larger example for parts
of the theory exposed in the Lectures on the dynamical Yang-Baxter equations
by P. Etingof and O. Schiffmann, math.QA/9908064. We explicitly compute
the matrix entries of the intertwining operator, fusion matrix and
exchange matrix associated to $sl(2,R)$ for $q=1$. We also consider
the weighted trace function and the dual Macdonald-Ruijsenaars equation
for this particular case. As expected, the matrix entries of the
exchange matrix turn out to be Racah polynomials. However, the road
to their explicit formula via the fusion matrix is quick, and it also
yields an interesting way to derive their orthogonality relations.

\Sec{intro} {Introduction}
The quantum dynamical Yang-Baxter equation (QDYBE)
was first considered in 1984
by Gervais and Neveu \ref{G-N}, with motivation from physics
(for monodromy matrices in Liouville theory).
A general form of QDYBE with spectral parameter was presented by
Felder \ref{Fe1}, \ref{Fe2} at two major congresses in 1994.
The corresponding classical dynamical Yang-Baxter equation (CDYBE)
was presented there as well.
Next Etingof and Varchenko started a program to give geometric interpretations
of solutions of CDYBE (see \ref{E-V1}) and of QDYBE (see \ref{E-V2})
in the case without spectral parameter. In the context of this
program they pointed out a method to obtain solutions of QDYBE by
the so-called exchange construction (see \ref{E-V3}). This uses,
for any simple Lie algebra $\gog$,
representation theory of $U(\gog)$ or of its quantized version $U_q(\gog)$
in order to define successively the intertwining operator, the fusion matrix
and the exchange matrix. The matrix elements of the intertwining operator
and of the exchange matrix generalize respectively the Clebsch-Gordan
coefficients and the Racah coefficients to the case where the
first tensor factor is a Verma module rather than a finite dimensional
irreducible module. The exchange matrix is shown to satisfy QDYBE.
Etingof and Varchenko also started in \ref{E-V4}
a related program to connect the above
objects with weighted trace functions and with solutions of the
(q-)Knizhnik-Zamolodchikov-Bernard equation (KZB or qKZB).
A nice introduction to the topics indicated above was recently given by
Etingof and O. Schiffmann \ref{E-S}.

The present paper is intended as a somewhat larger example for parts of
the theory in \ref{E-S}. In particular, while
Example 2 of Section 2 of \ref{E-S} explicitly computes the matrix elements
of the intertwining operator, the fusion matrix and the exchange matrix
in the case of the 2-dimensional
irreducible representation of $sl(2)$, $q=1$, we will do this
for irreps of $sl(2)$ of arbitrary finite dimension.
This is the topic of Sections 3--5, after some preliminaries in Section 2.
An immediate consequence of
the explicit form of the fusion matrix is an explicit expression of
the universal fusion matrix for $sl(2)$, $q=1$, see Section 6, and of the
related universal operator $Q(\la)$, see Section 7. In
Section 8 of \ref{E-S} this particular univeral fusion matrix is
obtained as a consequence of the ABBR equation in the $q=1$ case.
We conclude in Section 8 with the example of the weighted trace function
for $sl(2)$, $q=1$ and the dual Macdonald-Ruijsenaars equation satisfied
by it. This illustrates part of Section 9.2 of \ref{E-S}.

The aim of this paper is mainly tutorial, as a complement to
\ref{E-S} for those who like to see larger worked-out examples.
The explicit expressions we obtain for the matrix elements of the
intertwining operator and the fusion matrix are terminating
${}_3F_2(1)$'s and terminating balanced ${}_4F_3(1)$'s, as should be
the case since these matrix elements, rational in $\la\in\CC$, must
be analytic continuations of Clebsch-Gordan coefficients respectively
Racah coefficients for finite dimensional irreps of $sl(2)$.
Also, the weighted trace function turns out to be a Gegenbauer function
of the second kind, which analytically continues the Gegenbauer polynomial
occurring as a weighted trace function associated with finite dimensional
irreps of $sl(2)$.

We would like to point to one observation of particular interest.
The fusion matrix $J_{\de,\ga}(\la)$, explicitly obtained in Section 4,
is triangular and has matrix entries given by elementary expressions,
and the same is true for the inverse fusion matrix.
The exchange matrix $R_{\de,\ga}(\la)$
is defined in terms of the fusion matrix and a transpose of its inverse
by  equation (5.1). Thus it is a product
of an upper triangular and a lower triangular matrix, both with elementary
matrix entries. Therefore we arrive very quickly at a single sum
expression for the matrix entries of the exchange matrix, which can
moreover immediately be recognized as Racah polynomials. This approach
seems quicker than other approaches in literature to explicit expressions
of Racah coefficients. Moreover, the inverse of the exchange matrix can
also be explicitly computed from the above factorization in an easy way,
which yields biorthogonality relations for the matrix entries of the
exchange matrix. After some transformations this yields the
known orthogonality relations for the Racah polynomials involved.

\mLP
{\sl Acknowledgement}\quad
The first author was inspired by Pavel Etingof's lectures on the dynamical
Yang-Baxter equations at the London Mathematical Society Symposium
on Quantum groups in Durham, UK, July 1999.

\mLP
{\sl Notation}\quad We will write $a\vee b$ for $\max(a,b)$ and
$a\wedge b$ for $\min(a,b)$. For notation and results on hypergeometric
series we refer to \ref{A-A-R} and \ref{Erd1}.

\Sec {prelim} {Preliminaries}
Let $\gog$ be a simple complex Lie algebra and let $\goh\i\gog$ be a
Cartan subalgebra. Let $\De\i\goh^*$ be the associated root system
and let $\De^+$ be a choice of positive roots.
Define a partial order $\le$ on $\goh^*$ by putting $\la\le\mu$
if $\mu-\la$ is in the $\Zplus$ span of $\De^+$.
If $V$ is any $\gog$-module and $\la\in\goh^*$, denote by $V[\la]$
the weight space $\{v\in V\mid h\dt v=\la(h)\,v
\hbox{ for all $h\in\goh$}\}$.
For $\la\in\goh^*$ let $M_\la$ be the Verma module with highest
weight $\la$ and corresponding highest weight vector $x_\la$.
Recall that the tensor product $W\ten V$ of two $\gog$-modules
$W$ and $V$ becomes a $\gog$-module such that
$g\dt(w\ten v)=(g\dt w)\ten v+w\ten(g\dt v)$\quad ($g\in\gog$, $w\in W$,
$v\in V$). Then $W\ten V$ also becomes a $U(\gog)$-module and we have
$$
g^n\dt(w\ten v)=\sum_{j=0}^n {n\choose j}(g^j\dt w)\ten(g^{n-j}\dt v)\quad
(g\in\gog,\;w\in W,\;v\in V,\;n\in\Zplus).
\eqno\eq{4}
$$

Take $\gog:=sl(2,\CC)$ with basis
$$
e:=\pmatrix{0&1\cr0&0},\quad
h:=\pmatrix{1&0\cr0&-1},\quad
f:=\pmatrix{0&0\cr1&0}.
$$
Let $\goh$ be spanned by $h$. Then the map
$\la\mapsto\la(h)\colon \goh^*\to \CC$ identifies $\goh^*$ with $\CC$.
We have $[h,e]=2e$, $[h,f]=-2f$, $[e,f]=h$. Thus $\De=\{2,-2\}$.
Choose $2\in\De$ as the positive root.
In $\goh^*$ we have that $\la\le\mu$ iff $\mu-\la$ is an even nonnegative
integer. By induction with respect to $n$ we have the following identity
in $U(\gog)$:
$$
ef^n=f^ne+nf^{n-1}(h-n+1)\qquad(n\in\ZZ_{\ge1}).
\eqno\eq{35}
$$
For $\la\in\goh^*$ a basis of the Verma module $M_\la$ is given by the
elements $f^k\dt x_\la$ ($k\in\Zplus$). Note that
$f^k\dt x_\la$  has weight $\la-2k$. Clearly we have $e\dt x_\la=0$.
Let $V$ be a finite dimensional irreducible $sl(2,\CC)$-module and
let $0\ne v\in V[\be]$. Then
$V=\union_{k\in\ZZ}V[\be+2k]$ and there are $k_0,k_1\in\Zplus$
such that $\dim V[\be+2k]=1$ for
$k=-k_0,\ldots,k_1$ and
$\dim V[\be+2k]=0$ otherwise. If $k\in\Zplus$ then $V[\be+2k]$ is spanned
by $e^k\dt v$.

\Sec {twine} {The intertwining operator}
Let $V$ be a finite dimensional $\gog$-module and
let $v\in V$ be a weight vector. For $\la\in\goh^*$
put $\mu:=\la-\wt(v)$.
As shown in \reff{E-S}{Proposition 2.2}, there is
for generic $\la\in\goh^*$ a unique $\gog$-intertwining operator
$\Phi_\la^v\colon M_\la\to M_\mu\ten V$ such that
$$
\Phi_\la^v(x_\la)\in x_\mu\ten v+\bigoplus_{\nu<\mu}M_\mu[\nu]\ten V.
$$
Moreover, the operator $\Phi_\la^v$ depends rationally on $\la$.

\Lemma{1}
Let $V$ be a finite dimensional irreducible
$sl(2,\CC)$-module, $v\in V$, $\be:=\wt(v)$. Then
$$
\Phi_\la^v(x_\la)=\sum_{k=0}^\iy
{1\over k!\,(-\la+\be)_k}\,(f^k\dt x_{\la-\be})\ten (e^k\dt v)
\qquad(\la\in\CC\backslash(\be+\Zplus)).
\eqno\eq{2}
$$

\Proof
By \reff{E-S}{Proposition 2.2} and by our earlier remarks about
$sl(2,\CC)$-modules, we must have for generic $\la\in\CC$ that
$$
\Phi_\la^v(x_\la)=\sum_{k=0}^\iy a_k\,(f^k\dt x_{\la-\be})\ten (e^k\dt v)
$$
for certain coefficients $a_k$ depending on $\la$ and $\be$ with
$a_0=1$.
By the intertwining property of $\Phi_\la^v$ we have
$e\dt \Phi_\la^v(x_\la)=0$. Hence, by \eqtag{35},
$$
\eqalignno{
0&=\sum_{k=1}^\iy a_k\,(ef^k\dt x_{\la-\be})\ten (e^k\dt v)
+\sum_{k=0}^\iy a_k\,(f^k\dt x_{\la-\be})\ten (e^{k+1}\dt v)
\cr
&=\sum_{k=1}^\iy a_k\,k(\la-\be-k+1)\,(f^{k-1}\dt x_{\la-\be})\ten(e^k\dt v)
+\sum_{k=0}^\iy a_k\,(f^k\dt x_{\la-\be})\ten(e^{k+1}\dt v)
\cr
&=\sum_{k=1}^\iy(a_{k-1}-k(k-\la+\be-1)a_k)\,(f^{k-1}\dt x_{\la-\be})\ten
(e^k\dt v).
}
$$
Hence, 
for $\la\in\CC\backslash(\be+\Zplus)$ and
for $k_1$ such that $e^k\dt v\ne0$ iff $k\le k_1$, we obtain the recurrence
$$
a_k={a_{k-1}\over k(k-\la+\be-1)}\qquad
(k=1,\ldots,k_1).
$$
By iteration we find that
$a_k=\displaystyle{1\over k!\,(-\la+\be)_k}$.\hhalmos

\bPP
Let $V_\ga$ be the finite dimensional irreducible $sl(2,\CC)$-module with
highest weight $\ga\in\Zplus$. Then we can take a basis of $V_\ga$
consisting of vectors $v_\ga^\ga,v_{\ga-2}^\ga,\ldots,v_{-\ga}^\ga$
such that
$$
\eqalignno{
h\dt v_{-\ga+2k}^\ga&=(-\ga+2k)\,v_{-\ga+2k}^\ga\,,
\cr
f\dt v_{-\ga+2k}^\ga&=k\,v_{-\ga+2k-2}\qquad
\hbox{(or 0 if $k=0$),}
\cr
e\dt v_{-\ga+2k}^\ga&=(\ga-k)\,v_{-\ga+2k+2}\qquad
\hbox{(or 0 if $k=\ga$).}
}
$$
Then
$$
e^i\dt v_{-\ga+2k}^\ga=(-1)^i\,(-\ga+k)_i\,v_{-\ga+2k+2i}^\ga\,,\quad
f^i\dt v_{-\ga+2k}^\ga={k!\over (k-i)!}\,v_{-\ga+2k-2i}^\ga\,.
\eqno\eq{5}
$$
Thus formula \eqtag{2} can be rewritten as
$$
\Phi_\la^{v_{-\ga+2k}^\ga}(x_\la)=
\sum_{i=0}^{\ga-k}{(-1)^i\,(-\ga+k)_i\over i!\,(-\la-\ga+2k)_i}\,
(f^i\dt x_{\la+\ga-2k})\ten v_{-\ga+2k+2i}^\ga\,,
\eqno\eq{3}
$$
where $\ga\in\Zplus$, $k\in\{0,1,\ldots,\ga\}$,
$\la\in\CC\backslash\{-\ga+2k,-\ga+2k+1,\ldots,k\}$.

\bPP
Let $n\in\Zplus$ and apply $f^n$ to both sides of \eqtag{3}.
By the intertwining property of $\Phi$, by equation \eqtag{4} and by the
second part of equation \eqtag{5} we obtain
$$
\eqalignno{
\Phi_\la^{v_{-\ga+2k}^\ga}(f^n\dt x_\la)
=&\sum_{i=0}^{\ga-k}\sum_{j=0\vee (n-k-i)}^n
{n\choose j}
{(-1)^i\,(-\ga+k)_i\over i!\,(-\la-\ga+2k)_i}\,
{(k+i)!\over (k+i+j-n)!}\,\cr
&\qquad\qquad\qquad\qquad\qquad\qquad\times
(f^{i+j}\dt x_{\la+\ga-2k})\ten v_{-\ga+2(k+i+j-n)}^\ga
\cr
=&\sum_{m=0\vee(n-k)}^{n+\ga-k} c_{m,n}^{\la,\ga,k}\,(f^m\dt x_{\la+\ga-2k})
\ten v_{-\ga+2(k+m-n)}^\ga,&\eq{6}
}
$$
where
$$
\eqalignno{
c_{m,n}^{\la,\ga,k}&=
{n!\,k!\over (k+m-n)!}\,
\sum_{i=0\vee (m-n)}^{(\ga-k)\wedge m}
{(-1)^i\,(-\ga+k)_i\,(k+1)_i\over (-\la-\ga+2k)_i\,(n-m+i)!\,(m-i)!\,i!}
&\eq{7}
\cr
&={n!\,k!\over (k+m-n)!}\,
\sum_{j=0\vee(n-m)}^{(n-m+\ga-k)\wedge n}
{(-1)^{j+m-n}\,(-\ga+k)_{j+m-n}\,(k+1)_{j+m-n}\over
(-\la-\ga+2k)_{j+m-n}\,(m-n+j)!\,(n-j)!\,j!}\,.\qquad&\eq{10}
}
$$
We obtained the last part of \eqtag{6} (with \eqtag{7} substituted)
by passing from summation variables $i,j$ to summation variables
$m,i$ with $m=i+j$. Furthermore, \eqtag{10} is obtained from
\eqtag{7} by substituting $i=j+m-n$ for the summation variable.

In order to write \eqtag{7} or \eqtag{10} in terms of hypergeometric functions,
we distinguish two cases.

\Theor{36}
Let the coefficients of the intertwining operator be given by
$$
\Phi_\la^{v_{-\ga+2k}^\ga}(f^n\dt x_\la)=
\sum_{m=0\vee(n-k)}^{n+\ga-k} c_{m,n}^{\la,\ga,k}\,(f^m\dt x_{\la+\ga-2k})
\ten v_{-\ga+2(k+m-n)}^\ga.
\eqno\eq{37}
$$
Then
$$
\eqalignno{
&c_{m,n}^{\la,\ga,k}=
{n!\,k!\over m!\,(n-m)!\,(k+m-n)!}\,
{}_3F_2\left[{-m,-\ga+k,k+1\atop -\la-\ga+2k,n-m+1};1\right]\qquad\qquad
\cr
\noalign{\medskip}
&\qquad\qquad\qquad\qquad\qquad\qquad\qquad\qquad\qquad\qquad
{\rm if}\quad 0\vee(n-k)\le m\le n;\qquad&\eq{8}
\cr
\noalign{\bigskip}
&c_{m,n}^{\la,\ga,k}=
{(-1)^{m-n}\,(-\ga+k)_{m-n}\over (m-n)!\,(-\la-\ga+2k)_{m-n}}\,
{}_3F_2\left[{-n,-\ga+k+m-n,k+m-n+1\atop -\la-\ga+2k+m-n,m-n+1};1\right].
\cr
\noalign{\medskip}
&\qquad\qquad\qquad\qquad\qquad\qquad\qquad\qquad\qquad\qquad
{\rm if}\quad 0\le n\le m\le n+\ga-k.\qquad&\eq{9}
}
$$

We will later need the special case $n=0$ of \eqtag{9}
and the special case $m=0$ of \eqtag{8}
(which are both already an immediate consequence of \eqtag{3}):
$$
\eqalignno{
c_{m,0}^{\la,\ga,k}&={(-1)^m\,(-\ga+k)_m\over m!\,(-\la-\ga+2k)_m}\quad
(m=0,1,\ldots,\ga-k),&\eq{12}
\cr
c_{0,n}^{\la,\ga,k}&={k!\over (k-n)!}\quad(n=0,1,\ldots,k).&\eq{13}
}
$$

\Sec {fusion} {The fusion matrix}
We define the fusion matrix following \reff{E-S}{\S2.1}.
For $\gog$ a simple complex Lie algebra
let $V$ and $W$ be finite-dimensional $\gog$-modules, and
let $v\in V$ and $w\in W$
be weight vectors. Then, for generic $\la\in\goh^*$ there is a unique
weight vector $J_{WV}(\la)(w\ten v)\in W\ten V$ of weight $\wt (v)+\wt(w)$
such that
$$
(\Phi_{\la-\wt(v)}^w\ten 1)(\Phi_\la^v(x_\la))=
\Phi_\la^{J_{WV}(\la)(w\ten v)}(x_\la).
\eqno\eq{11}
$$
For generic $\la$ there is a unique extension of $J_{WV}(\la)$
to a $\goh$-linear map $J_{WV}(\la)\colon W\ten V\to W\ten V$.
The operator $J_{WV}(\la)$ is called the {\sl fusion matrix} of $V$ and $W$.
It depends rationally on $\la$.

We will compute the fusion matrix for $\gog=sl(2,\CC)$.
Let $\ga,\de\in\Zplus$, $k\in\{0,1,\ldots,\ga\}$,
$l\in\{0,1,\ldots,\de\}$. Let $\la\in\CC$ be generic.
Then, by \eqtag{37},
$$
\eqalignno{
&(\Phi_{\la+\ga-2k}^{v_{-\de+2l}^\de}\ten 1)
(\Phi_\la^{v_{-\ga+2k}^\ga}(x_\la))=
\sum_{r=0}^{\ga-k} c_{r,0}^{\la,\ga,k}\,\Phi_{\la+\ga-2k}^{v_{-\de+2l}^\de}
(f^r\dt x_{\la+\ga-2k})\ten v_{-\ga+2k+2r}^\ga
\cr
&\quad=\sum_{r=0}^{\ga-k}\sum_{m=0\vee(r-l)}^{r-l+\de}
c_{r,0}^{\la,\ga,k}\,c_{m,r}^{\la+\ga-2k,\de,l}\,
(f^m\dt x_{\la+\ga+\de-2k-2l})\ten v_{-\de+2(l+m-r)}^\de\ten
v_{-\ga+2k+2r}^\ga
\cr
&\quad=\sum_{m=0\vee(-l)}^{\ga-k-l+\de}f^m\dt x_{\la+\ga+\de-2k-2l}\ten
\sum_{r=0\vee(m+l-\de)}^{(\ga-k)\wedge(m+l)}
c_{r,0}^{\la,\ga,k}\,c_{m,r}^{\la+\ga-2k,\de,l}\,
v_{-\de+2(l+m-r)}^\de\ten v_{-\ga+2k+2r}^\ga\,.
}
$$
Write $J_{\de,\ga}$ instead of $J_{V_\de,V_\ga}$.
It follows by combination with \eqtag{11} and by substitution of
\eqtag{12}, \eqtag{13} that
$$
\eqalignno{
J_{\de,\ga}(\la)(v_{-\de+2\l}^\de\ten v_{-\ga+2k}^\ga)&=
\sum_{r=0}^{(\ga-k)\wedge l}
c_{r,0}^{\la,\ga,k}\,c_{0,r}^{\la+\ga-2k,\de,l}\,
v_{-\de+2l-2r}^\de\ten v_{-\ga+2k+2r}^\ga\,.
\cr
&=\sum_{r=0}^{(\ga-k)\wedge l}
{(-1)^r\,(-\ga+k)_r\,l!\over r!\,(-\la-\ga+2k)_r\,(l-r)!}\,
v_{-\de+2l-2r}^\de\ten v_{-\ga+2k+2r}^\ga.
}
$$
Replace $l,k,r$ by $n,s-n,n-m$, respectively.
Then we get for $0\le s\le\ga+\de$ and $(s-\ga)\vee 0\le n\le\de\wedge s$ that
$$
J_{\de,\ga}(\la)(v_{-\de+2n}^\de\ten v_{-\ga+2s-2n}^\ga)=
\sum_{m=(s-\ga)\vee 0}^n
A_{m,n}^{\la-\ga,\ga-s}\,
v_{-\de+2m}^\de\ten v_{-\ga+2s-2m}^\ga
\eqno\eq{14}
$$
where
$$
A_{m,n}^{\la-\ga,\ga-s}=A_{m,n}:=
{(-1)^{n-m}\,n!\,(-\ga+s-n)_{n-m}\over(n-m)!\,m!\,(-\la-\ga+2s-2n)_{n-m}}
\quad\hbox{if $m\le n$,}
\eqno\eq{15}
$$
and $A_{m,n}:=0$ if $m>n$.
Note that $A_{m,n}$ is independent of $\de$, and that, as a function
of $\la$, $\ga$, $s$, it depends only on $\la-\ga$ and $\ga-s$.
Note also that $A_{m,n}$ remains well-defined for generic complex values
of $\ga$ and $s$ (not necessarily integer)
and for $m,n\in\Zplus$. Thus $(A_{m,n})_{m,n\in\Zplus}$
is an infinite upper triangular matrix from which the matrix occurring
in \eqtag{14} is obtained by restricting $m$ and $n$ to the 
finite set $\{(s-\ga)\vee 0,\ldots,\de\wedge s\}$.

We will need later an explicit expression for the inverse $(B_{m,n})$
of the matrix $(A_{m,n})$.
By using Maple we found the explicit value of $(B_{m,n})$
for $m,n$ restricted to $\{0,\ldots,s\}$ with $s$ small.
From this we conjectured the general expression and next proved it:

\Lemma{16}
For $\la,\ga,s\in\CC$ generic and for $m,n\in\Zplus$
(and moreover $m,n\ge s-\ga$ if $s-\ga\in\ZZ_{>0}$)
let
$$
B_{m,n}^{\la-\ga,\ga-s}=B_{m,n}:=
{n!\,(-\ga+s-n)_{n-m}\over (n-m)!\,m!\,(-\la-\ga+2s-m-n-1)_{n-m}}
\quad\hbox{if $m\le n$,}
\eqno\eq{17}
$$
and $B_{m,n}:=0$ if $m>n$. Then the matrix $(B_{m,n})$
is the inverse of the matrix $(A_{m,n})$.

\Proof
For $m\le n$ we have
$$
\eqalignno{
&\sum_{k=m}^n B_{m,k}\,A_{k,n}=
\sum_{l=0}^{n-m} B_{m,m+l}\,A_{m+l,n}
\cr
&=\sum_{l=0}^{n-m}{(m+l)!\,(-\ga+s-m-l)_l\,(-1)^{n-m-l}\,n!\,
(-\ga+s-n)_{n-m-l}
\over
l!\,m!\,(-\la-\ga+2s-2m-l-1)_l\,(n-m-l)!\,(m+l)!\,(-\la-\ga+2s-2n)_{n-m-l}}
\cr
&\qquad
={(-1)^{n-m}\,n!\,(-\ga+s-n)_{n-m}\over m!\,(n-m)!\,(-\la-\ga+2s-2n)_{n-m}}
\,\sum_{l=0}^{n-m}{(-n+m)_l\,(\la+\ga-2s+n+m+1)_l\over
(\la+\ga-2s+2m+2)_l\,l!}
\cr
&\qquad
={(-1)^{n-m}\,n!\,(-\ga+s-n)_{n-m}\over m!\,(n-m)!\,(-\la-\ga+2s-2n)_{n-m}}
\,{}_2F_1\left[{-n+m,\la+\ga-2s+n+m+1\atop \la+\ga-2s+2m+2};1\right]
\cr
&\qquad\qquad\qquad
={(-1)^{n-m}\,n!\,(-\ga+s-n)_{n-m}\over m!\,(n-m)!\,(-\la-\ga+2s-2n)_{n-m}}
\,{(-n+m+1)_{n-m}\over(\la+\ga-2s+2m+2)_{n-m}}\,,
}
$$
which equals 0 if $m<n$ and equals 1 if $m=n$.
In the last identity of the displayed formula we used the
{\sl Chu-Vandermonde identity}
$$
{}_2F_1\left[{-n,b\atop c};1\right]={(c-b)_n \over (c)_n}\qquad (n\in\Zplus),
\eqno\eq{44}
$$
see for instance
\reff{A-A-R}{Corollary 2.2.3}. In particular, we have
$$
{}_2F_1\left[{-n,c+n-1\atop c};z\right]={(-n+1)_n \over (c)_n}=\de_{n,0}
\qquad (n\in\Zplus)
\eqno\eq{45}\quad\halmos
$$

\Remark{43}
Because of Lemma \thtag{16}, $A_{m,n}$ and $B_{m,n}$ as given by \eqtag{15}
and \eqtag{17} satisfy
$$
\sum_{k=m}^n A_{m,k}\,B_{k,n}=\de_{m,n}=\sum_{k=m}^n B_{m,k}\,A_{k,n}.
\eqno\eq{48}
$$
We showed the second identity in the proof of Lemma \thtag{16}.
We may verify the first identity independently as follows.
$$
\eqalignno{
&\sum_{k=m}^n A_{m,k}\,B_{k,n}=\sum_{l=0}^{n-m} A_{m,n-l}\,B_{n-l,n}=
{(-1)^{n-m}\,n!\,(-\ga+s-n)_{n-m}\over m!\,(n-m)!\,(-\la-\ga+2s-2n)_{n-m}}
\cr
&\qquad\qquad\qquad\quad\times
{}_3F_2\left[{-n+m,-\la-\ga+2s-2n-1,\thalf(-\la-\ga+2s-2n-1)+1\atop
-\la-\ga+2s-n-m,\thalf(-\la-\ga+2s-2n-1)};1\right].
}
$$
Now observe that, for $n\in\Zplus$, an application of \eqtag{44} yields:
$$
\eqalignno{
{}_3F_2\left[{ -n,b,\thalf b+1\atop b+n+1,\thalf b};1\right]
&={}_2F_1\left[{-n,b \atop b+n+1};1\right]
-{2n\over b+n+1}\,{}2F_1\left[{-n+1,b+1\atop b+n+2};1\right]
\cr
&={(n+1)_n\over (b+n+1)_n}\,-\,{2n\over b+n+1}\,
{(n+1)_{n-1}\over (b+n+2)_{n-1}}=\de_{n,0}.&\eq{46}
}
$$

\Remark{38}
The explicit matrix inversion in Lemma \thtag{16} is a special case
of Gould \reff{Gou}{Theorem 2}. This can be seen if we rewrite
\eqtag{15} as
$$
A_{k,n}={(\la+\ga-2s+1)_k\over(\ga-s+1)_k}\,
{n!\,(\ga-s+1)_n\over(\la+\ga-2s+1)_{2n}}\,
(-1)^k\,{n\choose k}\,
{-\la-\ga+2s-1-k\choose n},
$$
and \eqtag{17} as
$$
\eqalignno{
B_{k,n}=&{(\la+\ga-2s+1)_{2k}\over k!\,(\ga-s+1)_k}\,
{(\ga-s+1)_n\over(\la+\ga-2s+1)_n}\,
{-\la-\ga+2s-1-n\choose n}^{-1}
\cr
&\qquad\qquad\times(-1)^k{-\la-\ga+2s-1-2k\over -\la-\ga+2s-1-n-k}\,
{-\la-\ga+2s-1-n-k\choose n-k}.}
$$
See further generalizations of formulas for explicit matrix inverses
in Riordan \reff{Rio}{Ch.\ 2,3} and
Krattenthaler \ref{Kra}.

\bPP
We summarize our results in the following theorem.

\Theor{37}
Let $(s-\ga)\vee 0\le n\le\de\wedge s$. Let $A_{m,n}^{\la-\ga,\ga-s}$
and $B_{m,n}^{\la-\ga,\ga-s}$ be given by \eqtag{15} respectively \eqtag{17}
if $0\le m\le n$ and put them equal to 0 for other $m,n\in\Zplus$.
Then the fusion matrix for $sl(2)$ and its inverse are given by
$$
\eqalignno{
J_{\de,\ga}(\la)(v_{-\de+2n}^\de\ten v_{-\ga+2s-2n}^\ga)&=
\sum_{m=(s-\ga)\vee 0}^n
A_{m,n}^{\la-\ga,\ga-s}\,
v_{-\de+2m}^\de\ten v_{-\ga+2s-2m}^\ga\,,&\eqtag{14}
\cr
J_{\de,\ga}(\la)^{-1}(v_{-\de+2n}^\de\ten v_{-\ga+2s-2n}^\ga)&=
\sum_{m=(s-\ga)\vee 0}^n
B_{m,n}^{\la-\ga,\ga-s}\,
v_{-\de+2m}^\de\ten v_{-\ga+2s-2m}^\ga\,.&\eq{18}
}
$$

\Sec{exchange} {The exchange matrix}
Let $\gog$ be a simple complex Lie algebra
and let $V$ and $W$ be finite-dimensional $\gog$-modules.
Again following \reff{E-S}{\S2.1} we define the {\sl exchange matrix}
in terms of the fusion matrix by
$$
R_{VW}(\la):=J_{VW}(\la)^{-1}\,J_{WV}^{21}(\la)\qquad
\hbox{($\la\in\goh^*$ generic).}
\eqno\eq{19}
$$
Here $J^{21}:=PJP$ with $P(x\ten y):=y\ten x$.
Then $R_{VW}(\la)\colon V\ten W\to V\ten W$ is an $\goh$-intertwining
linear operator, rationally depending on $\la$.
It follows immediately from \eqtag{19} that
$$
R_{VW}(\la)^{-1}=(R_{WV}(\la))^{21}
\eqno\eq{23}
$$

We will compute the exchange matrix for $\gog=sl(2,\CC)$.
Let $\ga,\de\in\Zplus$, $k\in\{0,1,\ldots,\ga\}$,
$l\in\{0,1,\ldots,\de\}$. Let $\la\in\CC$ be generic.
Write $R_{\de,\ga}$ instead of $R_{V_\de,V_\ga}$.
It follows from \eqtag{14} that
$$
\eqalignno{
J_{\ga,\de}^{21}(\la)\,(v_{-\de+2s-2n}^\de\ten v_{-\ga+2n}^\ga)&=
P\,J_{\ga,\de}(\la)\,(v_{-\ga+2n}^\ga\ten v_{-\de+2s-2n}^\de)
\cr
&=\sum_{k=0\vee (s-\de)}^n A_{k,n}^{\la-\de,\de-s}\,
v_{-\de+2s-2k}^\de\ten v_{-\ga+2k}^\ga\,.&\eq{20}
}
$$
Hence, by \eqtag{19}, \eqtag{20} and \eqtag{18} we obtain
for $(s-\de)\vee 0\le n\le\ga\wedge s$ that
$$
\eqalignno{
R_{\de,\ga}(\la)(v_{-\de+2s-2n}^\de\ten v_{-\ga+2n}^\ga)&=
\sum_{k=0\vee (s-\de)}^n A_{k,n}^{\la-\de,\de-s}\,J_{\de,\ga}(\la)^{-1}\,
(v_{-\de+2s-2k}^\de\ten v_{-\ga+2k}^\ga)
\cr
&
=\sum_{k=0\vee(s-\de)}^n\,\sum_{m=k}^{\ga\wedge s}B_{s-m,s-k}^{\la-\ga,\ga-s}\,
A_{k,n}^{\la-\de,\de-s}\,v_{-\de+2s-2m}^\de\ten v_{-\ga+2m}^\ga
\cr
&
=\sum_{m=0\vee(s-\de)}^{\ga\wedge s} C_{m,n}^{\la,\ga,\de,s}\,
v_{-\de+2s-2m}^\de\ten v_{-\ga+2m}^\ga\,,
}
$$
where
$$
\eqalignno{
&C_{m,n}^{\la,\ga,\de,s}:=
\sum_{k=0\vee(s-\de)}^{m\wedge n} B_{s-m,s-k}^{\la-\ga,\ga-s}\,
A_{k,n}^{\la-\de,\de-s}&\eq{32}
\cr
&=
\sum_{k=0\vee(s-\de)}^{m\wedge n}
{(s-k)!\,(-\ga+k)_{m-k}\,(-1)^{n-k}\,n!\,(-\de+s-n)_{n-k}
\over (m-k)!\,(s-m)!\,(-\la-\ga+m+k-1)_{m-k}\,(n-k)!\,k!\,
(-\la-\de+2s-2n)_{n-k}}\,
\cr
&\qquad
={(-1)^m\,(-s)_m\,(-\ga)_m\over m!\,(-\la-\ga+m-1)_m\,(-\la-\de+2s-2n)_n}\,
\sum_{k=0\vee(s-\de)}^{m\wedge n}(\de-s+k+1)_{n-k}
\cr
&\qquad\qquad\qquad\qquad\qquad\times
{(-m)_k\,(-\la-\ga+m-1)_k\,(-n)_k\,(\la+\de-2s+n+1)_k
\over (-s)_k\,(-\ga)_k\,k!}\,.&\eq{22}
}
$$
Now we can write \eqtag{22} in terms of hypergeometric functions, where
we will distinguish two cases (for the second case replace the
summation variable $k$  in \eqtag{22} by $k+s-\de$).

\Theor{39}
Let the coefficients $C_{m,n}^{\la,\ga,\de,s}$
of the exchange matrix for $sl(2)$ be defined by
$$
\eqalignno{
&R_{\de,\ga}(\la)(v_{-\de+2s-2n}^\de\ten v_{-\ga+2n}^\ga)
=\sum_{m=0\vee(s-\de)}^{\ga\wedge s} C_{m,n}^{\la,\ga,\de,s}\,
v_{-\de+2s-2m}^\de\ten v_{-\ga+2m}^\ga\,.&\eq{21}
\cr
\noalign{\hbox{Then the following holds.}\smallskip}
&\hbox{If $s\le\de$ and $m,n\in \{0,\ldots,\ga\wedge s\}$ then}
\qquad\qquad\qquad\qquad\qquad\qquad\qquad\qquad\qquad\qquad\qquad
\cr
&\qquad\qquad C_{m,n}^{\la,\ga,\de,s}=
{(-1)^m\,(-s)_m\,(-\ga)_m\,(\de-s+1)_n
\over m!\,(-\la-\ga+m-1)_m\,(-\la-\de+2s-2n)_n}
\cr
&\qquad\qquad\qquad\qquad\qquad\times
{}_4F_3\left[{-m,-\la-\ga+m-1,-n,\la+\de-2s+n+1\atop
-s,-\ga,\de-s+1};1\right];&\eq{26}
\cr
\noalign{\medskip}
&\hbox{if $s\ge\de$ and $m,n\in \{s-\de,\ldots,\ga\wedge s\}$ then}
\qquad\qquad\qquad\qquad\qquad\qquad\qquad
\cr
&\qquad C_{m,n}^{\la,\ga,\de,s}=
{(-1)^{n+\de-s}\,(-\de)_{m+\de-s}\,(s-\ga-\de)_{m+\de-s}\,n!
\over (\la+\ga-2m+2)_{m+\de-s}\,(\la-s+n+1)_{n+\de-s}\,(m+\de-s)!\,(s-\de)!}
\cr
&\qquad\quad\times
{}_4F_3\left[{-m+s-\de,-\la-\ga-\de+s+m-1,-n+s-\de,\la-s+n+1\atop
-\de,-\ga-\de+s,s-\de+1};1\right].&\eq{24}
}
$$

\bPP
The ${}_4F_3$ in \eqtag{26} can be recognized as a {\sl Racah polynomial},
see \ref{W} and \reff{K-S}{\S1.2}. In the notation of
\reff{K-S}{(1.2.1)} we have
$$\eqalignno{
&{}_4F_3\left[{-m,-\la-\ga+m-1,-x,\la+\de-2s+x+1\atop
-s,-\ga,\de-s+1};1\right]
\qquad\qquad\qquad\qquad\qquad\qquad
\cr
&\qquad\qquad
=R_m(x(x+\la+\de-2s+1);-\ga-1,-\la-1,-s-1,\la+\de-s+1).
&\eq{31}
}
$$

From \eqtag{23} we obtain that $R_{\de,\ga}(\la)^{-1}=(R_{\ga,\de}(\la))^{21}$.
Hence it follows from \eqtag{21} that
$$
R_{\de,\ga}(\la)^{-1}\,(v_{-\de+2s-2n}^\de\ten v_{-\ga+2n}^\ga)=
\sum_{m=0\vee(s-\de)}^{\ga\wedge s} C_{s-m,s-n}^{\la,\de,\ga,s}\,
v_{-\de+2s-2m}^\de\ten v_{-\ga+2m}^\ga\,.
\eqno\eq{27}
$$

Let us restrict ourselves from now on, for convenience, to the case that
$s\le\ga\wedge\de$. Then we are in the situation of \eqtag{26} and we have
that $m,n\in\{0,1,\ldots,s\}$.
By \eqtag{26} and \eqtag{27} the matrix elements of
$R_{\de,\ga}(\la)^{-1}$ can be expressed in terms of
$$
\eqalignno{
C_{s-m,s-n}^{\la,\de,\ga,s}&=
{(-1)^{s-m}\,(-s)_{s-m}\,(-\de)_{s-m}\,(\ga-s+1)_{s-n}
\over (s-m)!\,(-\la-\de+s-m-1)_{s-m}\,(-\la-\ga+2n)_{s-n}}
\cr
&\qquad\times
{}_4F_3\left[{m-s,-\la-\de+s-m-1,n-s,\la+\ga-s-n+1\atop
-s,-\de,\ga-s+1};1\right].\quad
&\eq{25}
}
$$
It will turn out that the balanced ${}_4F_3$ in \eqtag{25} can be rewritten
in terms of the balanced ${}_4F_3$ in \eqtag{26} with $m$ and $n$
interchanged. This will follow by two successive applications
of Whipple's ${}_4F_3$ transform
$$
{}_4F_3\left[{-n,a,b,c\atop d,e,f};1\right]=
{(e-a)_n\,(f-a)_n\over (e)_n\,(f)_n}\,
{}_4F_3\left[{-n,a,d-b,d-c\atop d,1+a-e-n,1+a-f-n};1\right],
\eqno\eq{28}
$$
where $n\in\Zplus$ and $a+b+c-n+1=d+e+f$
(see for instance \reff{A-A-R}{Theorem 3.3.3}).
In fact, we have
$$
\eqalignno{
&{}_4F_3\left[{n-s,\la+\ga-s-n+1,m-s,-\la-\de+s-m-1\atop
-s,-\de,\ga-s+1};1\right]=
{(-\la+n)_{s-n}\over (-\de)_{s-n}}
\cr
&\times
{(-\la-\ga-\de+s+n-1)_{s-n}\over
(\ga-s+1)_{s-n}}\,
{}_4F_3\left[{n-s,\la+\ga-s-n+1,-m,\la+\de-2s+m+1\atop
-s,\la+\ga+\de-2s+2,\la-s+1};1\right]
}
$$
and
$$
\eqalignno{
&{}_4F_3\left[{-m,\la+\de-2s+m+1,n-s,\la+\ga-s-n+1\atop
-s,\la+\ga+\de-2s+2,\la-s+1};1\right]=
{(\ga-m+1)_m\over(\la+\ga+\de-2s+2)_m}
\cr
&\qquad\qquad\times
{(-\de+s-m)_m\over(\la-s+1)_m}\,
{}_4F_3\left[{-m,\la+\de-2s+m+1,-n,-\la-\ga+n-1
\atop -s,-\ga,\de-s+1};1\right].
}
$$
Hence
$$
\eqalignno{
&{}_4F_3\left[{n-s,\la+\ga-s-n+1,m-s,-\la-\de+s-m-1\,
\atop
-s,-\de,\ga-s+1};1\right]
\cr
&\qquad={(-\la-\ga-\de+s+n-1)_{s-n}\,(-\la+n)_{s-n}\,
(\ga-m+1)_m\,(-\de+s-m)_m
\over
(-\de)_{s-n}\,(\ga-s+1)_{s-n}\,(\la+\ga+\de-2s+2)_m\,(\la-s+1)_m}
\cr
&\qquad\qquad\qquad\qquad\times
{}_4F_3\left[{-m,\la+\de-2s+m+1,-n,-\la-\ga+n-1
\atop -s,-\ga,\de-s+1};1\right].
&\eq{29}
}
$$

From \eqtag{21} and \eqtag{27} we obtain
(for $s\le\ga\wedge\de$):
$$
\sum_{x=0\vee (s-\de)}^{\ga\wedge s}
C_{m,x}^{\la,\ga,\de,s}\,C_{s-x,s-n}^{\la,\de,\ga,s}=\de_{m,n}\,.
\eqno\eq{30}
$$
We can consider \eqtag{30} as biorthogonality relations between the two
systems of functions on $\{0,\ldots,s\}$ given by
$x\mapsto C_{n,x}^{\la,\ga,\de,s}$ and
$x\mapsto C_{s-x,s-n}^{\la,\de,\ga,s}$
($n\in\{0,\ldots,s\}$).
If we substitute \eqtag{26}
and \eqtag{25} in \eqtag{30}, next substitute \eqtag{29}, and
finally substitute \eqtag{31}
then, after a computation,
the resulting orthogonality relations for Racah polynomials
precisely coincide with those given in \reff{K-S}{(1.2.2)}
(replace $\al,\be,\ga,\de,N$ in \reff{K-S}{(1.2.2)} by
$-\ga-1,-\la-1,-s-1,\la+\de-s+1,s$).

\Sec {universal} {The universal fusion matrix}
From \eqtag{14} we can compute the universal fusion matrix for $sl(2)$.
For $\gog$ a simple complex Lie algebra the {\sl universal fusion matrix}
is a suitable generalized element $J(\la)$ ($\la\in\goh^*$)
of $U(\gog)\ten U(\gog)$ such that for all finite dimensional
$\gog$-modules $V,W$ we have $J(\la)\bigr|_{W\ten V}=J_{WV}(\la)$,
see \reff{E-S}{\S8.1}.

In fact, it follows from \eqtag{14}, \eqtag{5} and \eqtag{15} that
$$
\eqalignno{
&J_{\de,\ga}(\la)\,(v_{-\de+2n}^\de\ten v_{-\ga+2s-2n}^\ga)
\cr
&=
\sum_{m=(s-\ga)\vee0}^n A_{m,n}^{\la-\ga,\ga-s}\,
{(-1)^{n-m}\,m!\over n!\,(-\ga+s-n)_{n-m}}\,
f^{n-m}\dt v_{-\de+2n}^\de\;\ten\; e^{n-m}\dt v_{-\ga+2s-2n}^\ga
\cr
&=
\sum_{m=(s-\ga)\vee0}^n{1\over(n-m)!\,(-\la-\ga+2s-2n)_{n-m}}\,
(f^{n-m}\ten e^{n-m})\dt(v_{-\de+2n}^\de\ten v_{-\ga+2s-2n}^\ga).
}
$$
Hence
$$
\eqalignno{
J_{WV}(\la)\,(w\ten v)&=
\left(\sum_{n=0}^\iy
{1\over n!\,(-\la+\wt(v))_n}\,f^n\ten e^n\right)\dt (w\ten v)
\cr
&=\sum_{n=0}^\iy\left(f^n\ten{1\over n!\,(-\la+h-2n)_n}\,e^n\right)
\dt(w\ten v),
}
$$
if $w$ and $v$ are weight vectors in finite dimensional
$sl(2)$-modules $W$ resp.\ $V$.
\LP
Since
$(-\la+h-2n)_n=(-1)^n\,(\la-h+n+1)_n$ we obtain:

\Theor{40}
The universal fusion matrix $J(\la)$ for $sl(2)$ equals
$$
J(\la)=\sum_{n=0}^\iy f^n\ten{(-1)^n\over n!\,(\la-h+n+1)_n}\, e^n.
\eqno\eq{41}
$$

\bPP
Formula \eqtag{41} is in agreement with the formula at the end of \S8.1
in \ref{E-S}. Earlier, a quantum analogue of \eqtag{41}
was given by Babelon \ref{Bab},
see also Babelon, Bernard \& Billey \reff{B-B-B}{\S2}.

In a quite similar way we can derive from \eqtag{18} the universal
inverse fusion matrix $J(\la)^{-1}$. We successively obtain:
$$
\eqalignno{
J_{\de,\ga}(\la)^{-1}\,(v_{-\de+2n}^\de\ten v_{-\ga+2s-2n}^\ga)=
\sum_{m=(s-\ga)\vee0}^n{(-1)^{n-m}\over(n-m)!\,(-\la-\ga+2s-m-n-1)_{n-m}}\quad&
\cr
\times
(f^{n-m}\ten e^{n-m})\dt(v_{-\de+2n}^\de\ten v_{-\ga+2s-2n}^\ga);&
}
$$
$$
\eqalignno{
J_{WV}(\la)^{-1}\,(w\ten v)&=
\left(\sum_{n=0}^\iy
{(-1)^n\over n!\,(-\la+\wt(v)+n-1)_n}\,f^n\ten e^n\right)\dt (w\ten v)
\cr
&=\sum_{n=0}^\iy\left(f^n\ten{(-1)^n\over n!\,(-\la+h-n-1)_n}\,e^n\right)
\dt(w\ten v);
}
$$
so the universal inverse fusion matrix equals
$$
J(\la)^{-1}=\sum_{n=0}^\iy f^n\ten{1\over n!\,(\la-h+2)_n}\, e^n.
\eqno\eq{42}
$$
A quantum analogue of \eqtag{42} was given in \reff{B-B-B}{\S2}.

At least formally, it should hold now that
$$
J(\la)\,J(\la)^{-1}=1\ten 1=J(\la)^{-1}\,J(\la).
\eqno\eq{47}
$$
We can also verify these indentities independently, quite parallel to
the verification of the two identities in \eqtag{48}.
For the proof of the second identity in \eqtag{47} note that
$$
\eqalignno{
&J(\la)^{-1}\,J(\la)=
\sum_{k,l=0}^\iy f^{k+l}\ten
{1\over k!\,(\la-h+2)_k}\,e^k\,{(-1)^l\over l!\,(\la-h+l+1)_l}\,e^l
\cr
&=\sum_{n=0}^\iy\left(1\ten
\sum_{k=0}^n{(-1)^{n-k}\over k!\,(\la-h+2)_k\,(n-k)!\,(\la-h+n+k+1)_{n-k}}
\right)(f^n\ten e^n)=1\ten 1,
}
$$
since, by \eqtag{45},
$$
\eqalignno{
&\sum_{k=0}^n{(-1)^{n-k}\over k!\,(\la-h+2)_k\,(n-k)!\,(\la-h+n+k+1)_{n-k}}
\cr
&\qquad\qquad={(-1)^n\over n!\,(\la-h+n+1)_n}\,
{}_2F_1\left[{-n,\la-h+n+1\atop\la-h+2};1\right]=\de_{n,0}.
}
$$

For the proof of the first identity in \eqtag{47} note that
$$
\eqalignno{
&J(\la)\,J(\la)^{-1}=\sum_{k,l=0}^\iy f^{k+l}\ten
{(-1)^l\over l!\,(\la-h+l+1)_l}\,e^l\,{1\over k!\,(\la-h+2)_k}\,e^k
\cr
&=\sum_{n=0}^\iy\left(1\ten\sum_{l=0}^n{(-1)^l\over l!\,
(\la-h+l+1)_l\,(n-l)!\,(\la-h+2l+2)_{n-l}}
\right)(f^n\ten e^n)=1\ten 1,
}
$$
since by \eqtag{46},
$$
\eqalignno{
&\sum_{l=0}^n{(-1)^l\over l!\,(\la-h+l+1)_l\,(n-l)!\,(\la-h+2l+2)_{n-l}}
\cr
&\qquad\qquad={1\over n!\,(\la-h+2)_n}\,
{}_3F_2\left[{-n,\la-h+1,\thalf(\la-h+1)+1
\atop \la-h+n+2,\thalf(\la-h+1)};1\right]=\de_{n,0}.
}
$$

\Sec{Q} {The operator $Q(\lambda)$}
We follow \reff{E-S}{\S9.2} and \reff{E-V4}{\S1.2}
by
defining, for generic $\la\in\goh^*$,
a generalized element $Q(\la)$ in $U(\gog)$ in terms of the universal
fusion matrix as follows:
$$
Q(\la):=\bigl(m\circ P\circ(1\ten S^{-1})\bigr)\,J(\la).
\eqno\eq{49}
$$
Here $m$ denotes the multiplication operator $x\ten y\mapsto xy$ in
$U(\gog)$ and $S$ denotes the antipode in $U(\gog)$.
It follows from \eqtag{41} that $Q(\la)$ for $sl(2)$ takes the
following form:
$$
Q(\la)=\sum_{n=0}^\iy e^n\,{1\over n!\,(\la+h+n+1)_n}\,f^n.
\eqno\eq{50}
$$
For $\ga\in\Zplus$ let $Q_\ga(\la)$ denote $Q(\la)$ acting on $V_\ga$.
It follows from \eqtag{50} and \eqtag{5} that, for $k\in\{0,1,\ldots,\ga\}$,
$$
Q_\ga(\la)\, v_{-\ga+2k}^\ga=
\left(\sum_{n=0}^k{(-1)^n(-\ga+k-n)_n\,k!\over n!\,(\la-\ga+2k-n+1)_n\,(k-n)!}
\right)v_{-\ga+2k}^\ga
$$
The coefficient of $v_{-\ga+2k}^\ga$ on the right can be rewritten as
${}_2F_1(-k,\ga-k+1;-\la+\ga-2k;1)=(-\la-k-1)_k/(-\la+\ga-2k)_k$
by \eqtag{44}. So we have obtained:

\Theor{51}
The operator $Q(\la)$ acting on $V_\ga$ is explicitly given by
$$
Q_\ga(\la)\, v_{-\ga+2k}^\ga=
{(-\la-k-1)_k\over (-\la+\ga-2k)_k}\,v_{-\ga+2k}^\ga\quad
(k\in\{0,1,\ldots,\ga\}).
\eqno\eq{52}
$$

\Sec{trace} {Weighted trace functions}
Weighted trace functions for $q=1$ are defined in \reff{E-V4}{\S10.1}.
This is by analogy to or as a limit case of the $q$-case.
Weighted trace functions in the $q$-case are defined in
\reff{E-V4}{\S1.2} and in \reff{E-S}{\S9.2}. Let us state here once more
the definition of weighted trace function for $q=1$.

For $\la\in\goh^*$ and $U$ a $\gog$-module let $\exp_\la$ be the endomorphism
of $U$ sending a weight vector $u$ in $U$ to
$e^{\lan\la,\wt(u)\ran}\,u$.
Let $V$ be a finite dimensional $\gog$-module, let $B[0]$ be a basis
of $V[0]$ and let $v^*\in V^*[0]$ be the dual basis vector corresponding to
a basis vector $v$ in $B[0]$.
For $\mu\in\goh^*$ generic let
$$
\Phi_\mu^{V[0]}:=\sum_{v\in B[0]}\Phi_\mu^v\ten v^*
\eqno\eq{53}
$$
and define
$$
\Psi_V(\la,\mu):=\Tr\bigl|_{M_\mu}(\Phi_\mu^{V[0]}\circ
\exp_\la)\in V[0]\ten V^*[0].
\eqno\eq{54}
$$
Because we have taken the trace in \eqtag{54},
we will generally obtain an infinite
sum which is a priori a formal
power series in the variables $e^{-\lan\la,\al_i\ran}$, where
$\al_1,\ldots,\al_{\dim\goh}$ are the simple roots.

Let $\rho\in\goh^*$ be half the sum of the positive roots of $\gog$.
Let $\lan\;,\;\ran$ be the nondegenerate symmetric bilinear form
on $\goh^*$ induced, up to a constant factor, by the Killing form on $\gog$
such that $\lan\al,\al\ran=2$ if $\al$ is a long root.
Let the Weyl denominator be given by
$$
\de(\la):=e^{\lan\la,\rho\ran}\prod_{\al>0}(1-e^{-\lan\la,\al\ran}).
\eqno\eq{55}
$$
For $sl(2,\CC)$, where $\goh^*$ is identified with $\CC$ and 2 is the only
positive root, we get $\rho=1$, $\lan\la,\mu\ran=\thalf\la\mu$ and
$$
\de(\la)=e^{\half\la}-e^{-\half\la}.
\eqno\eq{59}
$$

Now define the weighted trace function by
$$
F_V(\la,\mu):=\bigl(\id_V\ten Q_{V^*}^{-1}(-\mu-\rho)\bigr)\,
\Psi_V(\la,-\mu-\rho)\,\de(\la).
\eqno\eq{56}
$$
This is again in $V[0]\ten V^*[0]$.

Let us compute $\Psi_\ga(\la,\mu)$, i.e.\ \eqtag{54}
for $\gog=sl(2,\CC)$ and $V=V_\ga$ with $\ga$ an even nonnegative integer.
From \eqtag{37} we obtain
$$
\Phi_\mu^{v_0^\ga}(\exp_\la\dt f^n\dt x_\mu)=
e^{\half\la(\mu-2n)}\sum_{m=0\vee(n-\half\ga)}^{n+\half\ga}
c_{m,n}^{\mu,\ga,\half\ga}\,(f^m\dt x_\mu)\ten v_{m-n}^\ga.
$$
Hence, by \eqtag{53},
$$
\Phi_\mu^{V_\ga[0]}(\exp_\la\dt f^n\dt x_\mu)=
\Phi_\mu^{v_0^\ga}(\exp_\la\dt f^n\dt x_\mu)\ten (v_0^\ga)^*.
$$
Hence, by \eqtag{54},
$$
\Psi_\ga(\la,\mu)=e^{\half\la\mu}
\sum_{n=0}^\iy c_{n,n}^{\mu,\ga,\half\ga}\,e^{-n\la},
\eqno\eq{57}
$$
where we omitted $v_0^\ga\ten(v_0^\ga)^*$ on the right-hand side,
since the one-dimensional vector space $V_\ga[0]\ten V_\ga^*[0]$ can
be identified with $\CC$. By \eqtag{8} we have
$$
\eqalignno{
\sum_{n=0}^\iy c_{n,n}^{\mu,\ga,\half\ga}\,e^{-n\la}&=
\sum_{n=0}^\iy e^{-n\la}\,
{}_3F_2\left[{-n,-\thalf\ga,\thalf\ga+1\atop -\mu,1};1\right]
\cr
&=\sum_{n=0}^\iy e^{-n\la}\sum_{k=0}^{n\wedge(\half\ga)}
{(-n)_k\,(-\half\ga)_k\,(\half\ga+1)_k\over(-\mu)_k\,k!\,k!}.
}
$$
Now interchange the two summations and next substitute $n=m+k$ for
the summation variable $n$. The above double sum becomes
$$
\sum_{k=0}^{\half\ga}
\left(\sum_{m=0}^\iy{(k+1)_m\over m!}\,(e^{-\la})^m\right)
{(-\half\ga)_k\,(\half\ga+1)_k\over(-\mu)_k\,k!}\,(-e^{-\la})^k.
$$
The inner sum converges absolutely for $\Re\la>0$. Therefore,
by dominated convergence, we see that for $\Re\la>0$
the interchange of summation
was justified and the sum in \eqtag{57} converges absolutely.
The inner sum equals $(1-e^{-\la})^{-k-1}$. So the double sum equals
$$
(1-e^{-\la})^{-1}\,\sum_{k=0}^{\half\ga}
{(-\half\ga)_k\,(\half\ga+1)_k\over(-\mu)_k\,k!}\,
(1-e^\la)^{-k}.
$$
So we obtain
$$
\eqalignno{
\Psi_\ga(\la,\mu)&=e^{\half\la\mu}\,(1-e^{-\la})^{-1}\,
{}_2F_1\left[{-\thalf\ga,\thalf\ga+1\atop -\mu};(1-e^\la)^{-1}\right]
&\eq{58}
\cr
&=e^{\half\la\mu}\,(1-e^{-\la})^{\half\ga}\,
{}_2F_1\left[{\thalf\ga-\mu,\thalf\ga+1;\atop -\mu};e^{-\la}\right],
&\eq{71}
}
$$
where we used Pfaff's transformation
${}_2F_1(a,b;c;z)=(1-z)^{-b}\,{}_2F_1\bigl(c-a,b;c;z/(z-1)\bigr)$,
see for instance \reff{A-A-R}{(2.2.6)}.
Combination with \eqtag{56}, \eqtag{59} and \eqtag{52} yields:

\Theor{60}
The weighted trace function for $sl(2,\CC)$ and $V_\ga$
($\ga$ nonnegative even integer) is given by
$$
\eqalignno{
F_\ga(\la,\mu)&=
{(-1)^{\half\ga}\,(\mu+1)_{\half\ga}\over(-\mu+1)_{\half\ga}}\,
e^{-\half\la\mu}\,
{}_2F_1\left[{-\thalf\ga,\thalf\ga+1\atop\mu+1};
(1-e^\la)^{-1}\right]
&\eq{61}
\cr
&={(-1)^{\half\ga}\,(\mu+1)_{\half\ga}\over(-\mu+1)_{\half\ga}}\,
e^{-\half\la\mu}\,(1-e^{-\la})^{\half\ga+1}\,
{}_2F_1\left[{\thalf\ga+\mu+1,\thalf\ga+1\atop\mu+1};
e^{-\la}\right].\qquad
&\eq{72}
}
$$

\bPP
Note that the restriction $\Re\la>0$ is no longer needed for convergence
in \eqtag{58} or \eqtag{61}.
The quantum analogue of \eqtag{61} was obtained in
\reff{E-V4}{Proposition 7.3}.

\mPP
By application to \eqtag{72} of the quadratic transformation formula
\reff{Erd1}{2.11(36)} (see also \reff{A-A-R}{(3.1.11)}),
and by next comparing with \reff{Koo84}{(2.15)},
we can express $F_\ga(\la,\mu)$ in terms of a {\sl Jacobi function
of the second kind}:
$$
\eqalignno{
&{(-1)^{\half\ga}\,(-\mu+1)_{\half\ga}\over(\mu+1)_{\half\ga}}\,
(e^{\half\la}-e^{-\half\la})^{-\half\ga-1}\,F_\ga(\la,\mu)
\cr
&\qquad\qquad=
(e^{{1\over 4}\la}+e^{-{1\over 4}\la})^{-\ga-2\mu-2}\,
{}_2F_1\left[{\thalf\ga+\mu+1,\mu+\thalf\atop 2\mu+1};
{1\over\cosh^2({1\over 4}\la)}\right]
\cr
&\qquad\qquad=\Phi_{2i\mu}^{(\half\ga+\half,\half\ga+\half)}(\tfrac 14 \la).
&\eq{73}
}
$$
\bPP
Finally we want to check \reff{E-S}{Theorem 9.2} (the dual
Macdonald-Ruijsenaars equations) for $sl(2,\CC)$ (so $q=1$).
The theorem states that for $V$ and $W$ finite-dimensional
$\gog$-modules we have
$$
\FSD_W^{\mu,V^*}\,F_V(\la,\mu)=\chi_W(e^{-\la})\,F_V(\la,\mu),
\eqno\eq{62}
$$
where
$$
\chi_W(e^\la):=\Tr\bigl|_W \exp_\la=\sum_\nu \dim(W[\nu])\,e^{\lan\la,\nu\ran}.
\eqno\eq{63}
$$
and $\FSD_W^{\mu,V^*}$
is a difference operator
$$
\eqalignno{
\FSD_W^{\mu,V^*}:=&\sum_{\nu\in\goh^*}\Tr|_{W[\nu]}\,(R_{WV^*}(-\mu-\rho))\,
T_\nu^\mu
\cr
=&\sum_{\nu\in\goh^*}\Tr|_{W[\nu]}\,
\bigl(R_{W[\nu],V^*[0];W[\nu],V^*[0]}(-\mu-\rho)\bigr)\,T_\nu^\mu.
&\eq{64}
}
$$
In \eqtag{64} $T_\nu^\mu$ denotes the shift operator defined by
$(T_\nu^\mu f)(\mu):=f(\mu+\nu)$.
Furthermore, in \eqtag{64}
$R_{W[\nu],V^*[\si];W[\nu'],V^*[\si']}$ denotes the block of the matrix
$\RR_{WV^*}$ corresponding to the weight spaces
$W[\nu],V^*[\si];W[\nu'],V^*[\si']$ (which block will be zero unless
$\nu+\si=\nu'+\si'$).

For the case of $sl(2,\CC)$ we obtain from \eqtag{64} and \eqtag{21} that
$$
\FSD_{V_\de}^{\mu,(V_\ga)^*}=
\sum_{s=\half\ga}^{\half\ga+\de}C_{\half\ga,\half\ga}^{-\mu-1,\ga,\de,s}\,
T_{-\de-\ga+2s}^\mu,
\eqno\eq{65}
$$
where $\ga$ is even and
the coefficients $C_{m,n}^{\la,\ga,\de,s}$ of the exchange matrix
are explicitly given by \eqtag{22}, or alternatively by
\eqtag{26} if $s\le\de$ or \eqtag{24} if $s\ge\de$.
Also, for $W=W_\de$, formula \eqtag{63} becomes
$$
\chi_\de(e^\la)=\sum_{k=0}^\de e^{\half\la(-\de+2k)}=
{e^{\half(\de+1)\la}-e^{-\half(\de+1)\la}\over e^{\half\la}-e^{-\half\la}}\,.
\eqno\eq{66}
$$
Thus the general theory yields that the difference equations
\eqtag{62} hold with \eqtag{61}, \eqtag{65} and \eqtag{66} substituted.
For general $\de$ we don't know if this formula was presented earlier.
However, for $\de=1$ we can reduce the formula to a well-known
contiguous relation for Gaussian hypergeometric series.

Indeed, for $\de=1$ we get
$$
C_{\half\ga,\half\ga}^{-\mu-1,\ga,1,\half\ga+1}=1;\quad
C_{\half\ga,\half\ga}^{-\mu-1,\ga,1,\half\ga}=
{(\mu-\half\ga-1)(\mu+\half\ga)\over(\mu-1)\mu}\,.
\eqno\eq{67}
$$
Thus
$$
F_\ga(\la,\mu+1)+{(\mu-\half\ga-1)(\mu+\half\ga)\over(\mu-1)\mu}\,
F_\ga(\la,\mu-1)=
(e^{\half\la}+e^{-\half\la})\,F_\ga(\la,\mu)
\eqno\eq{68}
$$
with $F_\ga(\la,\mu)$ given by \eqtag{61}.
This coincides (after appropriate substitution of parameters and argument)
with the contiguous relation
$$
\eqalignno{
&c(c-1)(z-1)\,{}_2F_1(a,b;c-1;z)+c\left(c-1-(2c-a-b-1)z\right)\,
{}_2F_1(a,b;c;z)
\cr
&\qquad\qquad\qquad\qquad\qquad\qquad\qquad
+(c-a)(c-b)\,{}_2F_1(a,b;c+1;z)=0,
&\eq{69}
}
$$
see \reff{Erd1}{2.8(30)}.

\Ref
\refitem A-A-R
G. E. Andrews, R. Askey \& R. Roy,
{\sl Special functions},
Cambridge University Press, 1999.

\refitem Bab
O. Babelon,
{\sl Universal exchange algebra for Bloch waves and Liouville theory},
Comm. Math. Phys. 139 (1991), 619--643.

\refitem B-B-B
O. Babelon, D. Bernard \& E. Billey,
{\sl A quasi-Hopf algebra interpretation of quantum 3-j and 6-j symbols
and difference equations},
Phys. Lett. B 375 (1996), 89--97;
{\tt \hbox{q-alg/9511019}}.

\refitem{Erd1}
A. Erd\'elyi e.a.,
{\sl Higher transcendental functions, Vol. I},
McGraw-Hill, 1953.

\refitem E-S
P. Etingof \& O. Schiffmann,
{\sl Lectures on the dynamical Yang-Baxter equations},
\item{}
preprint {\tt math.QA/9908064 v2}, 1999, 2000.

\refitem E-V1
P. Etingof \& A. Varchenko,
{\sl Geometry and classification of solutions of the classical dynamical
Yang-Baxter equation},
Comm. Math. Phys. 192 (1998), 77--120;
\item{}
{\tt q-alg/9703040}.

\refitem E-V2
P. Etingof \& A. Varchenko,
{\sl Solutions of the quantum dynamical Yang-Baxter equation and dynamical
quantum groups},
Comm. Math. Phys. 196 (1998), 591--640;
\item{} {\tt q-alg/9708015}.

\refitem E-V3
P. Etingof \& A. Varchenko,
{\sl Exchange dynamical quantum groups},
Comm. Math. Phys. 205 (1999), 19--52;
{\tt math.QA/9801135}.

\refitem E-V4
P. Etingof \& A. Varchenko,
{\sl Traces of intertwiners for quantum groups and difference equations, I},
{\tt math.QA/9907181}, 1999, 2000.

\refitem Fe1
G. Felder,
{\sl Conformal field theory and integrable systems associated to elliptic
curves},
in {\sl Proceedings of the International Congress of Mathematicians,
Z\"urich, 1994},
\item{}
Birkh\"auser 1994, pp. 1247--1255.

\refitem Fe2
G. Felder,
{\sl Elliptic quantum groups},
in
{\sl XIth International Congress of Mathematical Physics (Paris, 1994)},
Internat. Press, Cambridge, MA, 1995, pp. 211--218;
\item{} {\tt hep-th/9412207}.

\refitem G-N
J.-L. Gervais and A. Neveu,
{\sl Novel triangle relations and absence of tachyons in Liouville
string field theory},
Nucl. Phys. B 238 (1984), 125--141.

\refitem Gou
H. W. Gould,
{\sl A series transformation for finding convolution identities},
Duke Math.~J. 29 (1962), 393--404.

\refitem K-S
R. Koekoek and R.~F. Swarttouw,
{\sl The Askey-scheme of hypergeometric orthogonal polynomials and its
$q$-analogue},
Report 98-17, Faculty of Technical Mathematics and Informatics,
Delft University of Technology, 1998.

\refitem Koo84
T. H. Koornwinder,
{\sl Jacobi functions and analysis on noncompact semisimple Lie groups},
in {\sl Special functions: Group theoretical aspects and applications},
R. A. Askey, T. H. Koornwinder \& W. Schempp (eds.),
Reidel, 1984, pp. 1--85.

\refitem Kra
C. Krattenthaler,
{\sl A new matrix inverse},
Proc. Amer. Math. Soc. 124 (1996), 47--59.

\goodbreak
\refitem Rio
J. Riordan,
{\sl Combinatorial identities},
Wiley, 1968.

\refitem W
J. A. Wilson,
{\sl Some hypergeometric orthogonal polynomials},
SIAM J. Math. Anal. 11~(1980), 690--701.

\nonfrenchspacing
\vskip 0.5truecm
\parindent 1.0truein
\addressfont
T.~H. Koornwinder,
Korteweg-de Vries Institute, Universiteit van Amsterdam,

Plantage Muidergracht 24, 1018 TV Amsterdam,
The Netherlands;
\sPP
email: {\ttaddressfont thk@science.uva.nl}
\bPP
N. Touhami,
Korteweg-de Vries Institute, Universiteit van Amsterdam,

Plantage Muidergracht 24, 1018 TV Amsterdam,
The Netherlands;
\sPP
Laboratoire de Physique Th\'eorique, Universit\'e d'Es-S\'enia,
Oran,
Algeria;
\sPP
email: {\ttaddressfont touhami@science.uva.nl},\quad
{\ttaddressfont nabila\_$\,$touhami@yahoo.com}

\bye